%% file: ms.tex
\documentclass[12pt,a4paper]{article}
\usepackage{amsmath,amssymb,latexsym}
\usepackage[dutch,english]{babel}
\usepackage{enumitem}
\setlist{nosep}
\newtheorem{theorem}{Theorem}
\newtheorem{ass}{Assumption}

\newtheorem{lemma}{Lemma}

\input wi4202mac

\newcommand{\varsub}[2]{\ensuremath{{\rm Var}_{#2}\mspace{-1mu}\left ( #1 \right )}}

\newcommand{\cexpecsub}[3]{\ensuremath{{\mathbb E}_{#3}\mspace{-1mu}\left[#1 \, | \, #2 \right]}}
\newcommand{\Qsub}[2]{\ensuremath{{\mathbb Q}_{#2}\!\left( #1 \right)}}

\newcommand{\Mb}{\ensuremath{{M_\beta}}}
\newcommand{\numin}{\ensuremath{{\nu_{\min}}}}
\newcommand{\numax}{\ensuremath{{\nu_{\max}}}}
\newcommand{\mumin}{\ensuremath{{\mu_{\min}}}}
\newcommand{\muminstar}{\ensuremath{{\mu^\ast_{\min}}}}
\newcommand{\mustar}{\ensuremath{{\mu^{\ast}}}}
\newcommand{\mustarx}{\ensuremath{{\mu^{\ast}_x}}}
\newcommand{\mustary}{\ensuremath{{\mu^{\ast}_y}}}

\newcommand{\tauplus}{\ensuremath{{\tau^{+}}}}
\newcommand{\alphamax}{\ensuremath{{\alpha_{\max}}}}

\newcommand{\alphastar}{\ensuremath{{\alpha^{\ast}}}}
\newcommand{\mumax}{\ensuremath{{\mu_{\max}}}}
\newcommand{\piG}{\ensuremath{{\pi_{G}}}}

\newcommand{\Pb}{\ensuremath{{P_{\negmedspace\negthinspace\beta}}}}

\newcommand{\calB}{\ensuremath{{\cal{B}}}}
\newcommand{\calBp}{\ensuremath{{\cal{B}_+}}}

\newcommand{\qed}{\ensuremath{{\square}}}

\begin{document}

\title{Exponential convergence of adaptive importance sampling estimators for Markov chain expectations}
\author{Ludolf E.~Meester\\
Delft Institute of Applied Mathematics,\\
Delft University of Technology,\\
The Netherlands.}



\maketitle
\begin{abstract}
In this paper it is shown that adaptive importance sampling algorithms converge at exponential rate
for Markov chain expectation problems that admit a combination of a filtered estimator and a Markov zero-variance measure. 
It extends a chain of results---special purpose proofs were already known for several cases~\cite{KBCP99,baggerly2000exponential,desai2001adaptive}.
A recent paper~\cite{awad2013zero} provides a complete description of the class of combinations of Markov process expectations of path functionals and 
filtered estimators that admit zero-variance importance measures that retain the Markov property. In a way, this is the maximal class for which 
adaptive importance sampling algorithms might exhibit exponential convergence.
The main purpose of this paper is to prove that this is the case: for (most of) those combinations the natural adaptive importance sampling algorithm converges at exponential rate.
In addition, the applicability of general Markov chain theory for this purpose is discussed through the analysis of a counterexample presented in~\cite{desai2001simulation}.
\end{abstract}




\section{Introduction}
The chain of results on adaptive Monte Carlo algorithms exhibiting exponential convergence started in the field of particle transport modeling with the work of 
Booth~\cite{Booth85}, followed by that of Kollman et al~\cite{KBCP99} and Baggerly et al~\cite{baggerly2000exponential}. 
Particle behaviour is modeled by discrete time Markov chains, where so-called ``scores'' are incurred upon transitions.
The papers describe adaptive importance sampling algorithms that estimate the expected total score incurred before absorption and these algorithms are shown to converge 
exponentially fast, both in the case of discrete~\cite{KBCP99} and continuous~\cite{baggerly2000exponential} state spaces.

In addition to these algorithms, Desai and Glynn~\cite{desai2001simulation} discuss an eigenvalue problem which is solved using similar importance sampling techniques and also leads 
to exponential convergence (details to be found in Desai's thesis~\cite{desai2001adaptive}). Their main concern, however, is the apparent similarity between these algorithms contrasted by the special nature of their proofs.  It is felt that they should be amenable to general state space Markov chain theory and they report on some insight into the technical difficulties of that venture.

Awad et al~\cite{awad2013zero} perform a systematic investigation into the question ``For which combinations of Markov processes and (expectations of) path functionals do there
exist filtered estimators that admit zero-variance importance measures that are Markov as well?'' For the discrete time case, 
the answer is a class of (generalized) expected cumulative discounted rewards that extend the total-score functionals; this is described in detail in~\cite[\S 3]{awad2013zero}.
In some sense, this is the maximal class of problems for which one might look for exponentially converging adaptive importance sampling algorithms. 
The purpose of this paper is to provide a convergence proof for as many problems in this discrete time class as possible; 
an obvious limitation is the restriction to those functionals that can be evaluated in finite time, a restriction which is not necessary for the general existence 
question addressed in~\cite{awad2013zero}.

The structure of the paper is as follows. 
In Section~\ref{setup} the class of Markov process expectations is described, following the setup in~\cite[\S 3]{awad2013zero}, simplifying in one aspect and adding some assumptions, mostly similar to the ones familiar from~\cite{KBCP99,baggerly2000exponential}. 
Section~\ref{classQnu} describes the class of Markov importance sampling measures and the filtered estimators, 
some bounds on the likelihood ratio and on the termination time $\tau$ under these measures, as well as the zero-variance estimator and measure.
The adaptive importance sampling algorithm is stated in Section~\ref{AISalg}.
In Section~\ref{thetheorem} the exponential convergence theorem and its proof are given. Characteristic are two elements
that can be found in~\cite{KBCP99,baggerly2000exponential,desai2001adaptive} as well: proving that variance contraction occurs once the iteration process has entered 
a ``good'' set (a neighborhood of the zero-variance point) and showing that (a subset of) this good set is visited infinitely often with probability one. 
The contraction is proved by direct expansion of the filtered estimators at the zero-variance point, this in contrast with the other papers, where this is done by analysis on 
the system of recursive equations for the variances. For the proof of the visiting argument, it seems inevitable to consider two subsequent iteration steps, 
as the only known way is via a uniform version of the strong law of large numbers, applied to the collection of importance sampling estimators for a compact subset of parameters; 
this collection is proved to be uniformly integrable by a weak continuity argument. Using Markov's inequality, it is shown that the preceding step enters this compact subset with high probability. 

Auxiliary to the main proof, a theorem (Theorem~\ref{expconMC}) on the exponential convergence of general state space Markov chains is proved, 
a variant on a theorem in~\cite{desai2001adaptive}. 
A martingale argument shows that exponential convergence occurs if there is contraction of test function expectations on a good set and there exists a 
uniform positive lower bound on the probability of entering a specific subset of this good set from anywhere in the state space in two time steps.
This last theorem connects with the discussion in Section~\ref{MCperspective} on the application of general Markov chain theory for a general proof. 
It is argued that it is probably not sufficient that the kernel governing the iteration updating process is weakly continuous, given the uniform integrability that is required in the proof. This issue and some others are discussed in connection with a counterexample from~\cite{desai2001simulation}.
Section~\ref{eigenvalue} discusses the eigenvalue problem from~\cite{desai2001adaptive}. Strictly speaking, it falls outside the setting as described in Section~\ref{setup}, 
yet the developed tools suffice to construct a proof, which is presented here as there seem to be some problems with
the proof in~\cite{desai2001adaptive}.

\section{Markov-process expectations}
\label{setup}
Consider a general (Polish) state space $S$ and a transition kernel $P(x, \md y)$,
defining a Markov measure~\pr~and a Markov chain $X=(\infseq{X}{0}{1})$.
Let $K \subset S$ and define $\tau = \inf\{n\geq 0: X_n \in K \}$.
Let $s$ and $\beta$ be real functions on $S \times S$ such that
$s(x,y) \geq 0$ and $\beta(x,y) >0$; for convenience, assume $s(x,y)=0$ for $x \in K$.
Let 
\begin{equation}
\label{defY}
Y = \sum_{i=1}^\tau s(X_{i-1},X_i) B_i, \qquad \text{where $B_i=\prod_{j=1}^i \beta(X_{j-1},X_j)$,}
\end{equation}
and let $\mu: S \to \realR$ be defined by $\mu(x)=\expecsub{Y}{x}=\cexpec{Y}{X_0=x}$; 
note that $\mu(x)=0$ for $x \in K$. 

This is the setting as described in Section~3 of~Awad, Glynn, and Rubinstein~\cite{awad2013zero}, 
the following being the correspondences (right-hand sides are in the Awad et al notation): $Y=Z-f(X_0)$, $\mu(x)=u(x)-f(x)$ and $s(x,y)=g(x,y)+f(x)$.
In that paper it is demonstrated that $f$ and~$g$ may be replaced by (particular) $\tilde{f}$ and $\tilde{g}$, without changing the total reward $Z$, 
but leading to different (zero-variance) change of measure and estimator. 
In this paper $f$ and $g$ are considered given, thus fixing a specific combination of zero-variance change of measure and estimator. 
Here, the topic of interest is the convergence of an adaptive algorithm to the zero-variance solution, whence without loss of generality
the per-visit rewards can be absorbed in the per-transition rewards, simplifying the notation: the role $f$ and $g$ play is only through $s(x,y)=g(x,y)+f(x)$.

Suppose $\probsub{Y>0}{x}=0$ for some $x \in K^c$, then $Y=0$ ${\mathbb P}_x$-a.s.~and $\mu(x)=0$. Therefore,
without loss of generality, such~$x$ are assumed to be in $K$ so that $A = \{ x \in S: \mu(x)>0\}=K^c$.
Let~\calB~denote the set of bounded measurable functions on $S$, equipped with the sup-norm, and \calBp~its subset of nonnegative functions;
for $f \in \calB$ it is assumed that $f(x)=0$ for $x \in K$.
A transition operator is associated with the transition kernel $P$ through 
$P f(x) = \int_S f(y)P(x, \md y)$ for $f \in \calB$; for functions like $s$, for which $s(x, \cdot) \in \calB$ this definition is extended
as $P s(x) = \int_S s(x,y)P(x, \md y)$.
The $\beta$-weighted version of $P$ is defined and denoted by $P_\beta(x,\md y)= \beta(x,y) P(x, \md y)$.

Writing $Y_x$ for $Y$ with $X_0=x$, conditioning on $X_1$ yields the distributional equivalence 
\begin{equation}
\label{recurY}
Y_x \eqdist \beta(x,X_1) s(x,X_1) + \beta(x,X_1) Y_{X_1} \indicator{X_1 \in A}.
\end{equation}
The function $\mu$ is known to be the smallest nonnegative solution to the integral equation
\begin{equation}
\label{inteqmu}
u(x) = \int_S [s(x,y)+u(y)]\beta(x,y)  P(x, \md y), \quad x \in A,
\end{equation}
subject to the side condition $u(x)=0$ for $x \in K$.
Define $h(x)=\expecsub{\beta(x,X_1) s(x,X_1)}{x}=(\Pb s)(x)$, then the equation may be written as $u=h+\Pb u$, 
with formal solution $u = \sum_{n=0}^\infty \Pb^n h = \sum_{n=1}^\infty \Pb^n s$. 

The following boundedness assumptions are part of the sufficient conditions for the exponential convergence result.

\begin{ass}
\label{bddness}
\begin{enumerate}[label=\bf\alph*]
\item 
\label{supsfin}
$M_s = \sup_{x \in A,y \in S} s(x,y) < \infty$.
\item
\label{supbetafin}
$M_\beta = \sup_{x \in A,y \in S} \beta(x,y) < \infty$.
\item 
\label{infbetapos}
$m_\beta = \inf_{x \in A,y \in S} \beta(x,y) > 0$.
\item
\label{mufinite}
$\mumax = \sup_{x \in A} \mu(x) < \infty$.
\end{enumerate}
\end{ass}

The last assumption implies that $h$ and $\Pb \mu$ 
are bounded.

\begin{ass}
\label{Ptransient}
A finite $m$ exists such that $\gamma=\inf_{x \in A} \expecsub{s(X_{\tau-1},X_\tau); \tau \leq m}{x}>0$.
\end{ass}
The usual notation $\expec{X; C}$ for $\expec{X \mathbf{1}_{C}}$ is employed.
Since $s$ is bounded from above, Assumption~\ref{Ptransient} implies that $\pi = \inf_{x \in A}\probsub{\tau \leq m, 
s(X_{\tau-1},X_\tau)>0}{x}>0$, which means that $\tau$ is stochastically dominated by 
a geometric distribution with parameter $\pi_m=1-\sqrt[m]{1-\pi}$. Assumption~\ref{Ptransient} is sufficient to 
guarantee the termination of simulations of $Y$ in finite expected time, under $\pr$ but also under the importance 
sampling measures defined in the next section (proved in Lemma~\ref{tauQnuUB}).
An example with a three state Markov chain in~\cite[Example 1]{lecuyer2008approximate} where $\mu(i)$ is the 
expected number of transitions before absorption in $\Delta$ given $X_0=i$ illustrates that 
without the above condition $\tau=\infty$ might occur.
Assumption~\ref{bddness}\ref{infbetapos} implies $Y \geq \min(m_\beta^\tau,1) \,s(X_{\tau-1},X_\tau)$ for $X_0=x \in A$,
which in combination with Assumption~\ref{Ptransient} shows $\expecsub{Y}{x} \geq \min(m_\beta^m,1) \,\gamma>0$, 
proving


\begin{lemma}
\label{infmupos}
$\mumin=\inf_{x \in A} \mu(x) >0$.
\end{lemma}

Assumptions~\ref{bddness}\ref{supsfin} and \ref{Ptransient} correspond to Assumption 3.1 and 3.2 in~\cite{baggerly2000exponential};
\ref{bddness}\ref{supbetafin} and \ref{infbetapos} might perhaps be weakened at the expense of complicating the proof.
Lemma~\ref{infmupos} is similar to the first part of Proposition 3.1~\cite{baggerly2000exponential}, whereas its second part does not follow 
from the other assumptions if $M_\beta>1$ and is therefore included as Assumption~\ref{bddness}\ref{mufinite}.

\section{A class of Markov importance sampling measures}
\label{classQnu}
Let $\nu \in \calBp$. For $x \in A$ define
$g(x) = h(x)+(\Pb\nu)(x)$
and consider the kernel $Q_\nu$ defined by
\begin{equation}
\label{Qnu}
Q_\nu(x,\md y)=
\begin{cases}
\frac{[s(x,y)+\nu(y)] \beta(x,y)}{g(x)}P(x,\md y) & \text{for $x \in A$, $g(x)>0$},\\
P(x, \md y) & \text{otherwise}.
\end{cases}
\end{equation}
Define
\begin{equation}
\label{lnu}
l_\nu(x,y) = 
\begin{cases}
\frac{g(x)}{[s(x,y)+\nu(y)] \beta(x,y)} & \text{for $x \in A$, $g(x)>0$, $s(x,y)+\nu(y)>0$},\\
1 & \text{otherwise}.
\end{cases}
\end{equation}
Let
\begin{equation}
\label{defLn}
L^\nu_n = \prod_{i=1}^n l_\nu(X_{i-1},X_i), \qquad L_n^{-\nu}=(L_n^\nu)^{-1},
\end{equation}
and
\begin{equation}
\label{Ynu}
Y_\nu = \sum_{i=1}^\tau s(X_{i-1},X_i)\,B_i\,L^\nu_i.
\end{equation}
From the definition of $Q_\nu$ it is clear that if $x \in A$, $g(x)>0$, and $s(x,y)+\nu(y)=0$, then the transition
$x \to y$ is not sampled under the corresponding Markov chain measure ${\mathbb Q}_\nu$. However, under the following assumption,
this is of no consequence. The arguments below are based on Theorem~1 of~\cite{awad2013zero} and its proof.

\begin{ass}
\label{nupos}
$\nu(x)>0$ for $x \in A$.
\end{ass}

From \eqref{recurY} and \eqref{inteqmu} it can be seen that $\mu(x)>0$ implies $h(x)>0$ or $P(x,A)>0$, whence 
Assumption~\ref{nupos} implies $g(x)>0$ for $x \in A$ and definition~\eqref{Qnu} may be stated as
\begin{equation}
\label{Qnumod}
Q_\nu(x,\md y)=
\begin{cases}
\frac{[s(x,y)+\nu(y)] \beta(x,y)}{g(x)}P(x,\md y) & \text{for $x \in A$;}\\
P(x, \md y) & \text{otherwise}.
\end{cases}
\end{equation}
Under Assumption~\ref{nupos} (recall that $\nu$ is zero on K) $s(x,y)+\nu(y)=0$ occurs if and only if $y \in K$ and 
$s(x,y)=0$, i.e., on a (final) transition into~$K$ with zero reward; 
the transition makes no contribution to $Y$. 
Considering $n$-step paths,
to make a contribution to $Y$ at time $n$ the path must be in 
\begin{equation}
A_n(x)=\{(\seq{x}{1}{n}): x_i \in A, i <n; s(x_{n-1},x_n)>0\}
\end{equation}
and for paths in this set $0<L_n^\nu<\infty$.
Paths in 
\begin{equation}
\label{Nnx}
N_n(x)=\{(\seq{x}{1}{n}): x_i \in A, i <n; s(x_{n-1},x_n)=0\} 
\end{equation}
have a zero contribution at time $n$; these paths are not sampled under $\mathbb Q_\nu$.
The measures ${\mathbb P}_x$ and ${\mathbb Q}_{\nu,x}$ are equivalent
on $\{(\seq{X}{1}{n}) \in A_n(x)\} \cap {\cal F}_{n \wedge \tau}$ and
\begin{equation}
\begin{split}
\expecsub{s(X_{n-1},X_n) B_n \indicator{\tau \geq n}}{x} 
& = \expecsub{s(X_{n-1},X_n) B_n \indicator{\tau \geq n} \indicator{(\seq{X}{1}{n}) \in A_n(x)}}{x} \\
& = \expecsub{s(X_{n-1},X_n) B_n L_n^{\nu} \indicator{\tau \geq n} \indicator{(\seq{X}{1}{n}) \in A_n(x)}}{\nu, x} \\
& = \expecsub{s(X_{n-1},X_n) B_n L_n^{\nu} \indicator{\tau \geq n}}{\nu, x},
\end{split}
\end{equation}
where $\mathbb{E}_{\nu,x}$ denotes expectation under $\mathbb Q_\nu$ given $X_0=x$.
This relation is used in~\cite{awad2013zero} to show that Assumption~\ref{nupos} implies 
$\expecsub{Y_\nu}{\nu,x}=\mu(x)$ for $x \in A$.
Replacing definition~\eqref{lnu} by
\begin{equation}
\label{lnumod}
l_\nu(x,y) = 
\begin{cases}
\frac{g(x)}{[s(x,y)+\nu(y)] \beta(x,y)} & \text{for $x \in A$;}\\
1 & \text{otherwise},
\end{cases}
\end{equation}
one has $0<L_n^\nu<\infty$, ${\mathbb Q}_{\nu}$-a.s.; paths in $N_n(x)$ have $L_n=\infty$, but 
$\Qsub{L_n^\nu=\infty}{\nu,x}=0$.
From the above, the validity of the following is easily seen:
\begin{lemma}
\label{lnucont}
For $\nu$ satisfying Assumption~\ref{nupos}:
$L_n^\nu(\omega)$ is positive, finite, and continuous in $\nu$, for $\omega \not\in\{\tau \leq n, s(X_{\tau-1},X_\tau)=0\}$; 
the exceptional set is a ${\mathbb Q}_{\nu}$-null set; 
$L_n^{-\nu}(\omega)$ is finite and continuous in $\nu$, for every $n$ and $\omega$. 
\end{lemma}

%


\paragraph{The zero-variance measure.}
Since $\mu$ satisfies~\eqref{inteqmu}, one obtains from \eqref{Qnumod} and \eqref{lnumod} 
\begin{equation}
\label{ZVkernel}
Q_\mu(x,\md y)=
\begin{cases}
\frac{[s(x,y)+\mu(y)]\beta(x,y)}{\mu(x)}P(x,\md y)
& \text{for $x \in A$},\\
P(x, \md y) & \text{otherwise},
\end{cases}
\end{equation}
and
\begin{equation}
\label{lmu}
l_\mu(x,y) = 
\begin{cases}
\frac{\mu(x)}{[s(x,y)+\mu(y)] \beta(x,y)} & \text{for $x \in A$,}\\
1 & \text{otherwise}.
\end{cases}
\end{equation}
Algebra shows that $Y_{\mu,x}$ evaluates to $\mu(x)$, $\mathbb Q_{\mu,x}$-a.s., whence
$\mathbb Q_\mu$ is called the zero-variance measure.

\subsection{Bounds}
\label{bounds}
Since $\mu$ is bounded away from zero and infinity on $A$, it is permitted to assume that this holds for~$\nu$ 
as well and only consider 
the set $B=\{\nu: \numin \leq \nu(x) \leq \numax,  x \in A\}$, where
$0<\numin \leq \mumin$ and $\mumax \leq \numax < \infty$.

Consider $\nu \in B$ and $x \in A$. The following lower bounds are easily obtained:
\begin{equation}
l_\nu^{-1}(x,y) \geq  
\begin{cases}
\frac{m_\beta}{(M_s+\nu_{\max})M_\beta} \cdot \nu_{\min} & \text{for $y \in A$,} \\
\frac{m_\beta}{(M_s+\nu_{\max})M_\beta} \cdot s(x,y) & \text{for $y \in K$.}
\end{cases}
\end{equation}
Therefore, a constant $\kappa \in (0,1)$ exists such that for $\nu \in B$:
\begin{equation}
\label{LBonLn}
L_n^{-\nu} \geq  
\begin{cases}
\kappa^n & \text{for $n < \tau$,} \\
\kappa^n \, s(X_{\tau-1},X_\tau) & \text{for $n \geq \tau$.}
\end{cases}
\end{equation}
Since $(\Pb \nu)(x) \geq (\Pb \mu)(x) \numin/\mumax$ and $\mu= h + \Pb \mu$:
\begin{equation*}
g(x) \geq h(x) + \frac{\numin}{\mumax} [\mu(x) -h(x)] \geq \frac{\numin \mumin}{\mumax} >0,
\end{equation*}
whence
\begin{equation}
 l_\nu^{-1}(x,y) 
\leq  \frac{(M_s+\numax)\Mb \mumax}{\numin \mumin}.
\end{equation}
Therefore, a finite constant $H$ exists such that for $\nu \in B$:
\begin{equation}
\label{UBonLn}
L_n^{-\nu} \leq  H^n.
\end{equation}

\begin{lemma}
\label{tauQnuUB}
The distribution of $\tau$ under $\mathbb{Q}_\nu$ is stochastically dominated by a single 
geometric distribution, for all $\nu \in B$.
\end{lemma}
\textbf{Proof.}
From~\eqref{Qnumod} it follows that $s(X_{\tau-1},X_\tau)>0$, $\mathbb Q_\nu$-a.s.~on $\{\tau<\infty\}$,
whence for $m$ as in Assumption~\ref{Ptransient}, using~\eqref{LBonLn}:
\begin{equation}
\label{LBonQnufortau}
\begin{split}
\Qsub{\tau \leq m}{\nu,x} 
& = \Qsub{\tau \leq m, s(X_{\tau-1},X_\tau)>0}{\nu,x} \\
& = \expecsub{L^{-\nu}_\tau; \tau \leq m, s(X_{\tau-1},X_\tau)>0 }{x} \\
& \geq \kappa^{m} \, \expecsub{s(X_{\tau-1},X_\tau) ; \tau \leq m}{x} \\
& \geq \kappa^{m} \gamma>0.
\end{split}
\end{equation}
Under $\mathbb{Q}_{\nu,x}$, $\tau$ is therefore dominated by a geometric random variable with parameter 
$\piG=1-\sqrt[m]{1-\kappa^m\gamma}$.~\qed



\section{Adaptive importance sampling}
\label{AISalg}
In the iterative algorithm below, since an estimate is needed for each $\mu(x)$, the chain and $Y_\nu$ are simulated for a number of different starting points 
$X_0=x$.
If $A$ is finite, a simulation can be done for each $x \in A$,
but if $A$ is infinite, a finite-dimensional model for $\mu$ and a finite set 
$D=\{\seq{x}{1}{d}\}$ of starting points are needed.
For example, a collection of linearly independent basis functions \seq{b}{0}{p} may be given such that 
\begin{equation}
\label{regmodel}
\mu(x;\alpha)= b_0(x)+  B(x) \alpha,
\end{equation}
where $B(x)=[b_1(x), \ldots, b_p(x)]$, $\alpha=[\seq{\alpha}{1}{p}]^T$, and $\mu(\cdot\,; \alpha_{\rm true})=\mu$.
Given a current estimate~$\nu$, 
$R$ replicates of the chain are simulated under $\mathbb{Q}_{\nu}$ starting from $X_0=x_i$, 
resulting in $Y_i^{(j)}, j=1, \ldots, R$; and this for each $1 \leq i \leq d$.
It~is assumed that the parameter vector $\alpha$ can be obtained from 
$\{Y_{i}^{(j)}, 1 \leq i \leq d, 1 \leq j \leq R\}$ by a suitable estimation procedure.
For simplicity's sake the estimate is assumed to be a function of the averages ${\bf y}=(\bar{Y}_i, 1\leq i \leq d$),
i.e., $\hat{\alpha}=\alpha({\bf y})$ for some function $\alpha: \realR^d \to \realR^p$.
Some regularity assumptions on the estimation procedure are needed:
\begin{ass}
\label{estimation}
\begin{enumerate}[label=\bf\alph*]
\item ${\bf y} \mapsto \mu(\cdot\,; \alpha({\bf y}))$ is Lipschitz continuous from $\realR^d_+$ to \calBp;
\item $\mu(\cdot\,; \alpha({\bf y}))=\mu$ for ${\bf y}=(\mu(x), x \in D)$;
\item $\mu(\cdot\,; \alpha({\bf y})) \geq \epsilon$, for some $\epsilon>0$.
\end{enumerate}
\end{ass}
Assumptions~\ref{estimation}\textbf{a}~and~\ref{estimation}\textbf{b}~hold, for example, when $B(x)$ is bounded and $\alpha$ is estimated by ordinary 
least squares. The third assumption is natural because $\mu$ is bounded away from zero on $A$ (Lemma~\ref{infmupos}).

These are the steps of the iterative algorithm:
\begin{enumerate}[label=\bf\arabic*.]
\item Choose an initial estimate $\mu^{(0)}$. Set $n=0$.
\item For each $x \in D$ simulate $R$ independent replications of $Y_{\mu^{(n)}}$ under $\mathbb{Q}_{\mu^{(n)}}$ with $X_0=x$.
\item Fit model \eqref{regmodel} to the simulated values to obtain $\hat\alpha$.
\item Set $\mu^{(n+1)}=\mu(\cdot\,;\hat\alpha)$.
\item Set $n \to n+1$ and repeat from step 2 until convergence conditions are met.
\end{enumerate}

\section{The exponential convergence theorem}
\label{thetheorem}
\begin{theorem} 
\label{expconvthm}
Under Assumptions 1--4 an $R_0$ and a $\theta>1$ exist
such that $\theta^n \, \vnorm{\mu^{(n)}-\mu} \to 0$ a.s.~for $R \geq R_0$.
\end{theorem}

The proof consists of three main parts: a general theorem on exponential convergence for Markov chains,
and proofs that two of its conditions are satisfied; the first is a contractive property of the 
variance of $\mu^{(n)}$ near the zero-variance point; the second that, starting from any point $\nu$, 
the estimate will be within this contractive set in two iteration steps, with fixed positive probability. 
These results are presented first, followed by the formal proof of Theorem~\ref{expconvthm}.

\subsection{Exponential convergence for general state space Markov chains}
\begin{theorem} 
\label{expconMC}
Let $X=\{X_n: n\geq0\}$ be any Markov chain, with state space $W$. Suppose there exist $F \subset G 
\subset W$, a real-valued function $g : W \to [0, \infty)$, 
and constants $c \in [0,1)$, ${\delta}>0$, finite~$b$ and $\tilde{b}$ satisfying 
\begin{enumerate}[label=\bf\Alph*]
\item\label{ass1} $\expecsub{g(X_1)}{x} \leq c \, g(x)$ for $x \in G$.
\item\label{ass2}  $g(x) \geq b$ for $x \in G^c$.
\item\label{ass3}  $g(x) \leq \tilde{b}<b $ for $x \in F$.
\item\label{ass4}  $\probsub{X_2 \in F}{x} \geq {\delta}$ for $x \in W$.
\end{enumerate}
Then $\theta^{n} {g}(X_n) \convas 0$ for some $\theta >1$.
\end{theorem}

This theorem is very similar to Theorem 5.3 in Desai~\cite[page 65]{desai2001adaptive}, with as crucial distinction the appearance of $X_2$ instead of $X_1$ in the last condition,
which prompts a number of changes in the proof, as \emph{two} time steps of the chain are involved rather than just one.

\paragraph{Proof.}
Define $\tilde{g}(x)=\min(g(x),b)$. Let $x \in G$. Assumption~\ref{ass1} implies
$c \, g(x)  \geq  \expecsub{g(X_1); X_1 \in G}{x} +  \expecsub{g(X_1); X_1 \not\in G}{x}$,
whence by assumption~\ref{ass2}:
\begin{equation}
\label{bndgX1}
\expecsub{g(X_1); X_1 \in G}{x} \leq c \, g(x) - b \, \probsub{X_1 \in G^c}{x}.
\end{equation}
Now
\begin{equation*}
\begin{split}
\expecsub{{g}(X_2); X_1 \in G}{x} 
& = \expecsub{\cexpecsub{g(X_2)\, \indicator{X_1 \in G}}{X_1}{x}}{x} \\
& = \expecsub{\indicator{X_1 \in G} \cexpecsub{g(X_2)}{X_1}{x}}{x} \\
& \leq \expecsub{\indicator{X_1 \in G} c \, g(X_1)}{x}  \\
& = c \, \expecsub{{g}(X_1); X_1 \in G}{x},
\end{split}
\end{equation*}
where the inequality follows from assumption~\ref{ass1}. Using $\tilde{g}(x) \leq b$ and 
the last two results:
\begin{equation*}
\begin{split}
\expecsub{\tilde{g}(X_2)}{x} 
& \leq c\, \expecsub{{g}(X_1); X_1 \in G}{x} + b \, \probsub{X_1 \in G^c}{x} \\
& \leq c^2 g(x) + (1-c)b \, \probsub{X_1 \in G^c}{x} \\
& = c g(x) \leq c \tilde{g}(x),
\end{split}
\end{equation*}
where $\probsub{X_1 \in G^c}{x} \leq {c \,g(x)}/{b}$ by Markov's inequality and assumption~\ref{ass1} and~\ref{ass2}.

Now let $x \in G^c$. Then 
\begin{equation*}
\begin{split}
\expecsub{\tilde{g}(X_2)}{x} 
& \leq \expecsub{{g}(X_2); X_2 \in F}{x} + b \, \probsub{X_2 \in F^c}{x}\\
&  \leq \tilde{b} \, \probsub{X_2 \in F}{x} + b \, \probsub{X_2 \in F^c}{x}\\
& = b - (b-\tilde{b})\, \probsub{X_2 \in F}{x} \\
& \leq b \, (1-(1-\tilde{b}/b){\delta}) 
= \tilde{g}(x) (1-(1-\tilde{b}/b){\delta}).
\end{split}
\end{equation*}
Set $\beta^2 = \max(c,1-(1-\tilde{b}/b){\delta})$, then $\expecsub{\tilde{g}(X_2)}{x} \leq \beta^2 \tilde{g}(x)$ for all $x$, 
with $\beta<1$, hence 
$\beta^{-2n}\tilde{g}(X_{2n})$ is a nonnegative supermartingale. By the martingale convergence 
theorem $\beta^{-2n} \tilde{g}(X_{2n}) \convas Z$ for some finite random variable~$Z$. 
Hence, 
$\tilde{g}(X_{2n}) \convas 0$ and $X_{2n} \in G$ eventually, i.e., also
$\beta^{-2n}{g}(X_{2n}) \convas Z$ and $\theta^{2n} {g}(X_{2n}) \convas 0$, for $1 < \theta < \beta^{-1}$.

It remains to show that this also holds for the odd indices. The above results imply that, eventually, $g(X_{2n}) \leq \beta^{2n} 
(2Z+1)$ 
and $\cexpec{g(X_{2n+1})}{X_{2n}} \leq c g(X_{2n})$. Hence, using Markov's inequality, for 
$\epsilon>0$:
\[
\prc{g(X_{2n+1})> \epsilon \theta^{-2n}}{X_{2n}} \leq \frac{c\, \theta^{2n}}{\epsilon}g(X_{2n}) \leq \frac{c}{\epsilon} 
(\theta \beta)^{2n} (2Z+1) \quad \text{eventually,}
\]
which implies
\[
 \sum_{n=1}^\infty \prc{g(X_{2n+1})> \epsilon \theta^{-2n}}{X_{2n}} < \infty \quad\text{a.s.}
\]
since $\theta\beta<1$. Then, by the conditional version of the Borel-Cantelli 
lemma~\cite[Cor.~5.29]{breiman1992probability},
$\prob{g(X_{2n+1})> \epsilon \theta^{-2n}~\text{i.o.}}=0$, i.e., $\theta^{2n} g(X_{2n+1}) \convas 0$. 
Combined with the above this proves $\theta^{n}g(X_n) \convas~0$.~\qed

\subsection{Error contraction near $\mu$}
\label{contraction}

\begin{lemma}
\label{contraclem}
A $\delta>0$ and a constant $M$ exist such that for $x \in A$
\begin{equation}
\label{varYnubound}
\varsub{Y_\nu}{\nu,x} \leq M \vnorm{\nu-\mu}^2 \qquad \text{whenever $\vnorm{\nu-\mu}<\delta$.}
\end{equation}
\end{lemma}
\textbf{Proof.}
Consider the difference $Y_{\nu,x} - Y_{\mu,x}$ under ${\mathbb Q}_{\nu,x}$. Under Assumption~\ref{nupos}, 
$L_n^\nu$ is positive and finite almost surely, for any $n$ and $\nu$ (Lemma~\ref{lnucont}). 
From \eqref{Ynu}:
\[
 Y_{\nu,x} - Y_{\mu,x}  =  
\sum_{n=1}^\tau s(X_{n-1},X_n) L_n^\mu B_n  \left (\frac{L_n^\nu}{L_n^\mu} -1 \right )
\]
and one obtains
\begin{equation}
\label{Ynu-mu}
 |Y_{\nu,x} - \mu(x)|   \leq  \mu(x) \cdot 
\max_{n \leq \tau} \left | \frac{L_n^\nu}{L_n^\mu} -1 \right |.
\end{equation}
A bound is sought on $L_n^\nu/L_n^\mu$ in terms of the error $\epsilon=\nu-\mu$. 
Define for $x \in A$
\begin{equation*}
\xi(x,y)= \frac{\epsilon(y)}{s(x,y)+\mu(y)},
\end{equation*}
setting $\xi(x,y)=0$ for $y \in K$.
Then
\begin{equation*}
\begin{split}
\epsilon(y)\Pb(x, \md y) 
& = \xi(x,y)[s(x,y)+\mu(y)]\Pb(x, \md y) \\
& =\mu(x) \xi(x,y) Q_\mu(x, \md y)
\end{split}
\end{equation*}
and, 
since $\mu$ satisfies \eqref{inteqmu},
\[
\int_{S} [s(x,z)+\nu(z)]\Pb(x, \md z)=
\mu(x)+(\Pb\epsilon)(x).
\]
Substituting the last two equations, one obtains:
\begin{equation}
\label{lnuratio}
\begin{split}
\frac{l_\nu(x,y)}{l_\mu(x,y)} 
& = \frac{\int_{S} [s(x,z)+\nu(z)] \Pb(x,\md z)}{[s(x,y)+\nu(y)] \beta(x,y)} \cdot \frac{[s(x,y)+\mu(y)]\beta(x,y)}{\mu(x)} \\
& = \frac{\mu(x)+ \Pb\epsilon(x)}{s(x,y)+\mu(y)+\epsilon(y)} \cdot \frac{s(x,y)+\mu(y)}{\mu(x)} \\
& = \Big ( 1 + \int_{S} \xi(x,z)Q_\mu(x, \md z)\Big ) \left ( 1 + \xi(x,y) \right ) ^{-1}.
\end{split}
\end{equation}
These manipulations are admissible, except when $x \in A$, $y \in K$, and $s(x,y)=0$; however, these combinations
are not sampled under ${\mathbb Q}_{\nu}$.
Since
\[
\me^{-\eta} \leq  \frac{1+a}{1+b}  \leq \me^\eta \qquad \text{for $|a|,|b| \leq  \eta/3$ and $\eta \leq 2/3$},
\]
one may conclude that
\begin{equation}
\label{boundlnu}
\me^{-\eta} \leq \frac{l_\nu(x,y)}{l_\mu(x,y)} \leq \me^\eta 
\end{equation}
provided $|\xi(x,y)| \leq \eta/3$ for all $x$, $y$, and $\eta \leq 2/3$.
Under this assumption 
\begin{equation}
\label{boundLnu}
\me^{-\eta n} \leq \frac{L_n^\nu}{L_n^\mu} \leq \me^{\eta n}
\quad\text{${\mathbb Q}_{\nu}$-a.s.}
\end{equation}
and then \eqref{Ynu-mu}
implies
\[
 |Y_{\nu,x} - \mu(x)|   
\leq \mu(x) \cdot (\me^{\eta\tau}-1)
\quad\text{${\mathbb Q}_{\nu}$-a.s.}
\]
By Lemma~\ref{tauQnuUB}, the distributions of $\tau$ under ${\mathbb Q}_{\nu,x}$ for $\nu \in B$ are uniformly bounded from above by
a geometric distribution with parameter $\piG$ (see the proof of Lemma~\ref{tauQnuUB}), which implies that for nonnegative $\eta$ the moment generating function
\expecsub{\me^{\eta\tau}}{\nu,x} is bounded above by the moment generating function $\phi(\eta)$ of that geometric distribution, therefore exists and is finite whenever $\me^\eta (1-\piG)<1$.
Thus,
\[
\expecsub{(\me^{\eta\tau}-1)^2}{\nu,x} \leq  \eta^2 \phi''(0)+o(\eta^2),
\]
when $\eta < -\log(1-\piG)$.
This implies, for $\eta$ small enough, that

\begin{equation}
\label{relvarbound}
\varsub{Y_\nu}{\nu,x} \leq \mu(x)^2 \, \expecsub{(\me^{\eta^\tau}-1)^2}{\nu,x}  \leq
2 \phi''(0) \mu(x)^2 \eta^2.
\end{equation}
Now, set $\eta= 3 \vnorm{\nu-\mu}/\mumin$, then $|\xi(x,y)| \leq |\nu(y)-\mu(y)| / \mu(y)$ implies $|\xi(x,y)| \leq \eta/3$ to imply~\eqref{boundlnu}, hence if \vnorm{\nu-\mu} is small enough, the last expression is bound by $M \vnorm{\nu-\mu}^2$ for some constant $M$.~\qed

\subsection{Close to the zero-variance point in two steps}

\begin{lemma}
\label{twosteps}
For every $\epsilon>0$ there exists and integer $R_0$ such that, if $R \geq R_0$ in the iterative algorithm in Section~\ref{AISalg}, then 
\[
 \prc{\vnorm{\mu^{(n+2)}-\mu} < \epsilon}{\mu^{(n)}} \geq \frac14 \quad\text{a.s.}
\]
\end{lemma}
Two lemma's are used in the proof:
\begin{lemma}
\label{weakcont}
$\nu \mapsto {\cal L}(\text{$Y_\nu$ under ${\mathbb Q}_\nu$})$ is continuous on $B$ in the topology of weak convergence.
\end{lemma}

\textbf{Proof.}
Let $\nu \in B$ and $f : \realR \to \realR$ continuous and bounded. The proof is complete if \expecsub{f(Y_\rho)}{\rho} is shown to be
continuous in $\rho$ at $\rho=\nu$.

By Lemma~\ref{lnucont}, 
all the terms of $\sum_{k=1}^\tau s(X_{k-1},X_k) \, B_k \, L^\rho_k$ are continuous, except perhaps the last one when  $s(X_{\tau-1},X_\tau)=0$, but in that case
the product $s(X_{\tau-1},X_\tau) \, L^\rho_\tau$ equals zero. Therefore, $Y_\rho(\omega)$ is continuous in $\rho$ at $\rho=\nu$, for all $\omega$.
Again by Lemma~\ref{lnucont}, $L_n^{-\rho}(\omega)$ is continuous in $\rho$ at $\rho=\nu$ for all $\omega$ and every $n$; by \eqref{UBonLn} it is
bounded from above by $H^n$.
Hence, $f(Y_\rho ) L_\tau^{-\rho} \indicator{\tau \leq n}$, 
is continuous in $\rho$ at $\rho=\nu$ and bounded uniformly in $\omega$, for any $n$,
so by the dominated convergence theorem 
\begin{equation}
\label{fromDCT}
\expec{f(Y_\rho ) L_\tau^{-\rho}; {\tau \leq n}}
\end{equation}
is continuous in $\rho$ at $\rho=\nu$, for all $n$. 
Writing $Z_\rho=f(Y_\rho ) L_\tau^{-\rho}$, one has
\[
\left | \expec{f(Y_\rho ) L_\tau^{-\rho}} - \expec{f(Y_\nu ) L_\tau^{-\nu}} \right | 
\leq 
\left | \expec{Z_\rho-Z_\nu; \tau \leq n} \right | +  \expec{|Z_\rho-Z_\nu|; \tau > n}
\]
and to prove the continuity, the second term needs to be controlled.
It is bounded by a constant multiple of 
$\expec{L_\tau^{-\rho} + L_\tau^{-\nu}; \tau >n}$.
For $\rho$ in a neighborhood of $\nu$, noting that $L_\tau^{-\rho}=0$ if $s(X_{\tau-1},X_\tau)=0$:
\[
\expecsub{L_\tau^{-\rho}; \tau \leq m}{x} = 
\expecsub{L_\tau^{-\rho}; \tau \leq m, s(X_{\tau-1},X_\tau)>0}{x} = Q_{\rho,x}(\tau \leq m)
\geq \kappa^m \gamma,
\]
as in the proof of Lemma~\ref{tauQnuUB},
for $m$ as in Assumption~\ref{Ptransient}.
Hence, $\expecsub{L_\tau^{-\rho}; \tau > n}{x} = Q_{\rho,x}(\tau > n) \to 0$ uniformly in $\rho$, by the same lemma.
Combined with the continuity of~\eqref{fromDCT} for fixed $n$, this implies that $\expec{f(Y_\rho ) L_\tau^{-\rho}}=\expecsub{f(Y_\rho)}{\rho}$
is continuous in~$\rho$ at $\rho=\nu$.~\qed

\begin{lemma}
\label{UIoverB}
The family $\{{\cal L}(\text{$Y_\nu$ under ${\mathbb Q}_\nu$}): \nu \in B\}$ is uniformly integrable.
\end{lemma}

\textbf{Proof.}
Let $i_n(y)$ be a continuous approximation of $\indicator{y \leq n}$, going linearly from one to zero on $[n,n+1]$, then $g_n(y)=y i_n(y)$ is bounded and continuous. 
By the previous lemma $\nu \mapsto \expecsub{g_n(Y_\nu)}{\nu,x}$ is continuous and since $\expecsub{Y_\nu}{\nu,x}=\mu(x)$ this
is also true for $\nu \mapsto f_n(\nu)=\expecsub{Y_\nu-g_n(Y_\nu)}{\nu,x}$.
However, $f_n(\nu) \to 0$ for $n \to \infty$, for any $\nu$ and $x$, so 
by Dini's theorem this implies that $f_n \to 0$ uniformly on the compact set $B$, i.e., $\expecsub{Y_\nu; Y_\nu>n}{\nu,x} \to 0$ uniformly on $B$.~\qed

\bigskip

\textbf{Proof of Lemma~\ref{twosteps}.} Following the steps that lead to~\cite[formula (42)]{baggerly2000exponential}, it can be shown from Markov's inequality that
$\prc{\mu^{(n+1)} \in B}{\mu^{(n)}} \geq \frac12$ a.s. Then, since 
$\{{\cal L}(\text{$Y_\nu$ under ${\mathbb Q}_\nu$}): \nu \in B\}$ is uniformly integrable, by Chung's uniform SLLN~\cite{chung1951strong}, a sample size $R_0$ exists such that
$\prc{\vnorm{\mu^{(n+2)}-\mu} < \epsilon}{\mu^{(n+1)}} \geq \frac12$ almost surely on the event $\{\mu^{(n+1)} \in B\}$, when $R \geq R_0$. Combining the two probilities leads to the required result.~\qed

\subsection{Proof of the convergence theorem}
\label{expconproof}
The proof is an application of Theorem~\ref{expconMC}. Clearly, $\{\mu^{(n)}: n \geq 0\}$ is a Markov chain, with state space \calBp.
Define $g(\nu) = \vnorm{\nu - \mu}^2$.
Suppose $\mu^{(0)}=\nu$ is given and $\bf y$ is the vector of averages $\bar{Y}_i$ from the simulations with $X_0=x_i$, for $x_i \in D$.
Then $\mu^{(1)}=\mu(\cdot\,;\alpha({\bf y}))$ and by Lipschitz continuity
\[
 \vnorm{\mu^{(1)}-\mu}^2 \leq C \sum_{i=1}^d \left (\bar{Y}_i-\mu(x_i) \right )^2
\]
for some constant $C$.
Hence
\begin{equation*}
\begin{split}
\cexpec{g\left ( \mu^{(n+1)} \right ) }{\mu^{(n)}} 
& = \cexpec{\vnorm{\mu^{(n+1)}-\mu}^2}{\mu^{(n)}} \\
& \leq C \sum_{i=1}^d \cexpec{\left (\bar{Y}_i-\mu(x_i) \right )^2}{\mu^{(n)}} \\
& \leq \frac{C}{R} \sum_{i=1}^d \varsub{Y_{\mu^{(n)}}}{\mu^{(n)},x_i}  \\
& \leq \frac{C d M}{R} \vnorm{\mu^{(n)}-\mu}^2 \\
& = c g\left ( \mu^{(n)} \right ),
\end{split}
\end{equation*}
whenever $\vnorm{\mu^{(n)}-\mu}<\delta$, with $\delta$ and $M$ as in~\eqref{varYnubound}. 
Furthermore, set
$G=\{\nu: g(\nu)<\epsilon\}$,
$F=\{\nu: g(\nu)\leq \epsilon/2\}$,
$R>C d M$ so $c<1$, then by choosing 
$b=\epsilon$ conditions \ref{ass1} through \ref{ass3} of Theorem~\ref{expconMC} are satisfied. Finally, choose $\epsilon$ as in Lemma~\ref{twosteps} and
condition~\ref{ass4} is also satisfied. The conclusion from Theorem~\ref{expconvthm} is the required result.~\qed


\section{An eigenvalue problem}
\label{eigenvalue}
In his PhD thesis~\cite{desai2001adaptive}, Desai analyzes adaptive importance sampling algorithms for finding the Perron-Frobenius eigenvalue of a 
finite nonnegative irreducible matrix $P$, assumed to be substochastic, without loss of generality. 
Labeling rows and columns by $0, 1, \ldots d$, the matrix $P$ defines a Markov chain on the state
space $S=\{0, 1, \ldots, d\}$, with the row defects assigned to transitions to the cemetary state $\Delta$. 
For the moment, define $K=\{\Delta\}$ and $\tau$ as usual.
State $0$ is chosen as \emph{return state} and $\sigma = \inf \{ n \geq 1:  X_n=0 \}$. For $\alpha\geq 0$ define
\begin{equation}
\label{Yalpha}
Y(\alpha)= \me^{\alpha \sigma} \indicator{\sigma<\tau}
\end{equation}
and $\mu_x(\alpha)=\expecsub{Y(\alpha)}{x}$.
The aim is to find $\alphastar$ satisfying $\mu_0( \alphastar)= \expecsub{Y(\alphastar)}{0}=1$, since then for $\mustar=\mu(\alphastar)$ and any $x$
\begin{equation}
\label{eigeneq}
 \mustarx= \expecsub{\cexpecsub{Y}{X_1}{x}}{x}=\me^{\alphastar} P(x,0)+\me^{\alphastar} \sum_{y=1}^d 
P(x,y)\mustary= \me^{\alphastar} \sum_{y=0}^d P(x,y)\mustary,
\end{equation}
that is, $\me^{-\alphastar}$ and $\mustar$ are the Perron-Frobenius eigenvalue and eigenvector. 

\medskip

An iteration of the adaptive algorithm consists of two updating steps. 
In the first, using the current estimate for $\mustar$, an importance sampling simulation estimate $\hat\mu_0(\alpha)$ is obtained for the function $\alpha \mapsto \mu_0(\alpha)$;
the next estimate for $\alphastar$ is the root $\hat\alpha$ of $\hat\mu_0(\alpha)=1$.
In the second step, updates for $\mustarx$ for $x \not \in \{ 0, \Delta\}$ are obtained by performing importance sampling simulations for~$\mu_x( \hat\alpha)$.
So, given an iterate $\mu^{(n)}$, in two updating steps the next iterate $\mu^{(n+1)}$ is obtained. Theorem 5.1~\cite{desai2001adaptive} reads (this paper's notation):

\begin{theorem}[Desai]
\label{thmdesai}
An $R_0$ and a $\theta>1$ exist such that $\theta^n \, \vnorm{\mu^{(n)}-\mustar} \to 0$ a.s.~for $R \geq R_0$.
\end{theorem}

As such, the algorithm does not fit into the framework described in this paper, because of the updating of $\alpha$. For fixed $\alpha$ it does fit in; this is described next.
In the setup of Section~\ref{setup}, define $s(x,y)=\delta_{y,0}$, $\beta(x,y)=\me^\alpha$, and $K=\{0,\Delta\}$. 
With these choices, for $X_0=x \not\in K$, \eqref{defY} is distributed as \eqref{Yalpha}. 
Denoting by $P_0$ the matrix $P$ after removal of row and column 0, by~\cite[Chapter 1, Theorem 5.1]{minc1988nonnegative} the largest eigenvalue of $P_0$ is strictly less than that of $P$.
This implies that an $\alphamax>\alphastar$ exists such that $\expecsub{Y(\alpha)}{x}<\infty$  for $\alpha \leq \alphamax$, 
so Assumption~\ref{bddness}\ref{mufinite} is satisfied for $\alpha\leq\alphamax$; parts \ref{supsfin}--\ref{infbetapos} are trivial. 
Assumption~\ref{Ptransient} is clearly satisfied for any $\alpha$, because $P$ is irreducible.
For $x \in A=K^c$, \eqref{Qnumod} and~\eqref{lnumod} become
\begin{equation}
\label{desaiQnu}
Q_\nu(x,y)= \frac{\delta_{y,0}+\nu(y)}{P(x,0)+\sum_{y=1}^d P(x,y)\nu(y)}P(x,y) \quad\text{and}\quad
l_\nu(x,y)= 
\frac{P(x,0)+\sum_{y=1}^d P(x,y)\nu(y)}{\delta_{y,0}+\nu(y)}.
\end{equation}
For fixed $\alpha$ the exponential convergence of $\mu^{(n)}(\alpha)$ to $\mu(\alpha)$ follows from 
Theorem~\ref{expconvthm}. 

\medskip

The tools developed, however, are sufficient to provide a convergence proof when $\alpha$ is updated.
The reason to present one is that there is a clear deficiency in the proof in~\cite{desai2001adaptive}:
at the bottom of page 72 it is asserted that (in current notation) $\expecsub{Y_{\nu}(\alpha)}{\nu,x}=\mu_x(\alphastar)$---i.e., 
with $\alphastar$ instead of $\alpha$---a claim that is incompatible with the strict increasingness of $Y_{\nu}(\alpha)$ as a function of $\alpha$ on the event $\{ Y_{\nu}(\alpha)>0\}$;
subsequently, the proof of Theorem 5.5~\cite{desai2001adaptive} uses the uniform integrability of the collection $\{\mu_x^{(n)}( \alpha^{(n)}) | \mu^{(n-1)}(\alpha^{(n-1)}) \in B \}$, which is undisputable, 
but it also apparently uses that $\cexpec{\mu_x^{(n)}( \alpha^{(n)})}{\mu^{(n-1)}(\alpha^{(n-1)})}=\mustarx$, probably based on the result on page 72. It is necessary, however, to show that the bias in
$\mu_x^{(n)}( \alpha^{(n)})$ vanishes sufficiently quickly, which is a pivotal issue.

\medskip

\textbf{Proof of Theorem~\ref{thmdesai}.}
The proof invokes Theorem~\ref{expconMC}, using $g(\nu)= \vnorm{\nu-\mustar}^2$, 
$G=\{\nu: g(\nu) < b\}$, $F=\{\nu: g(\nu) < b/2\}$.
The first and major part of the proof is showing that the contraction condition \ref{ass1} holds if $b$ is chosen small 
enough. Condition~\ref{ass4} is proved as in Lemma~\ref{twosteps} with some modifications. 

Updating~$\alpha$ involves $Y(\alpha)$ with $X_0=0$, which is zero because $0 \in K$ so $\tau=0$. To reproduce~\eqref{Yalpha} properly for this case 
define $Y_0(\alpha)=\me^{\alpha \sigma} \indicator{\sigma \leq\tau^+}$, where 
$\tau^+ = \inf \{ n \geq 1:  X_n \in K \}$ and modify $Y_{0,\nu}(\alpha)$ accordingly.
Note that the expressions in~\eqref{desaiQnu} do not depend on $\alpha$ as $\beta(x,y)$ is constant and cancels out, 
so neither do the likelihood ratios---only for $\nu=\mustar$ an~$\alphastar$ may appear through the application of~\eqref{eigeneq}. 
Under any $\mathbb{Q}_\nu$, $X_{\tau^+}=0$ and $\sigma=\tauplus$ a.s., whence the indicator can be dropped and
$Y_{0,\nu}(\alpha)= \me^{\alpha \tau^+} L_{\tau^+}^\nu$ a.s. (since the $\mathbb{Q}_\nu$ are equivalent measures, this 
is with respect to any of them).
Some algebra shows that $L^\mustar_\tauplus=\me^{-\alphastar\tauplus}$ and hence
\begin{equation}
\label{Ynualpha}
Y_{0,\nu}(\alpha)=
\frac{L_{\tau^+}^\nu}{L_{\tau^+}^\mustar} \,
\me^{(\alpha-\alphastar)\tauplus} \quad \text{a.s.}
\end{equation}

Repeating the steps in the proof 
of Lemma~\ref{contraclem} to bound $L^\nu_n/L^\mustar_n$ 
(replacing $\mu$ by $\mustar$ throughout) 
some changes occur in the derivation but the final result of~\eqref{lnuratio} is identical, whence also
\begin{equation}
\label{Ltauplusbnd}
\me^{-\eta\tauplus} \leq \frac{L_{\tau^+}^\nu}{L_{\tau^+}^\mustar} \leq
\me^{\eta\tauplus} \quad \text{a.s.}
\end{equation}
for $\eta=3 \vnorm{\nu-\mustar}/\muminstar$, provided \vnorm{\nu-\mustar} is small enough. Now, consider replications
$Y_{0,\nu}^{(i)}(\alpha)$, $i=1, \ldots, R$ and their average $\hat\mu_0( \alpha)$, and note that these are pathwise 
strictly increasing and strictly convex functions of $\alpha$, since $\tauplus \geq 1$.
This also implies that $\alpha \mapsto \mu_x(\alpha)$ is Lipschitz for all $x$, so a $\partial$ exists such that 
$\vnorm{\mu(\alpha)-\mu(\alphastar)} \leq \partial |\alpha -\alphastar|$. 
Combining the last two results shows that the root $\hat\alpha$ of $\hat\mu_0( \alpha)=1$ satisfies 
$|\hat\alpha-\alphastar|\leq \eta$. A bound on the second moment is also required, which is developed now. 
Noting that $\mu_0(\alphastar)=\hat\mu_0(\hat\alpha)=1$, $\hat\alpha$ satisfies
\begin{equation}
\mu_0(\alphastar)-\hat\mu_0(\alphastar)=
\hat\mu_0(\hat\alpha)-\hat\mu_0(\alphastar)=\hat\mu_0'(\xi)(\hat\alpha-\alphastar), 
\end{equation}
for some $\xi$ between $\hat\alpha$ and \alphastar, and
one obtains
\begin{equation}
\label{sqerr}
(\hat\alpha-\alphastar)^2 = \frac{[\hat\mu_0(\alphastar)-\mu_0(\alphastar)]^2}{[\hat\mu_0'(\xi)]^2}.
\end{equation}
From~\eqref{Ynualpha} and~\eqref{Ltauplusbnd} follow
\begin{equation}
[\hat\mu_0(\alphastar)-\mu_0(\alphastar)]^2 \leq \left ( \frac1R \sum_{i=1}^R \me^{\eta\tau^+_i}-1 \right )^2
\end{equation}
and
\begin{equation}
\hat\mu_0'(\xi) \geq 
\hat\mu_0'(\alphastar-\eta) =  
\frac1R \sum_{i=1}^R \tau^+_i Y_{0,\nu}^{(i)}(\alphastar-\eta) \geq
\frac1R \sum_{i=1}^R \me^{-2\eta\tau^+_i} \geq
 \left (\frac1R \sum_{i=1}^R \me^{2\eta\tau^+_i} \right )^{-1},
\end{equation}
where the last step follows from Jensen's inequality. Combining these with~\eqref{sqerr}:
\begin{equation}
 (\hat\alpha-\alphastar)^2  \leq \left ( \frac1R \sum_{i=1}^R \me^{\eta\tau^+_i}-1 \right )^2 \left (\frac1R \sum_{i=1}^R \me^{2\eta\tau^+_i} \right )^2.
\end{equation}
Writing $W_i=\me^{\eta\tau^+_i}$ and $V_i=\me^{2\eta\tau^+_i}$ and $\bar{W}$ and $\bar{V}$ for their averages, application of Cauchy-Schwarz implies
\begin{equation}
\expecsub{(\hat\alpha-\alphastar)^2}{\nu,0} \leq \sqrt{\expecsub{\bar{W}^4}{\nu,0} \expecsub{\bar{V}^4}{\nu,0}}.
\end{equation}
It is known that $\expecsub{\bar{W}^4}{\nu,0} = R^{-2} \left (\expecsub{(W-1)^2}{\nu,0} \right )^2 + o(R^{-2})$, provided $\expecsub{(W-1)^4}{\nu,0} < \infty$, 
and that \expecsub{\bar{V}^4}{\nu,0} converges to a constant as $R \to \infty$,  
provided the fourth moment of $V$ is finite. Thus,
\begin{equation}
\expecsub{(\hat\alpha-\alphastar)^2}{\nu,0} \leq \frac{c_1}{R} \expecsub{(W-1)^2}{\nu,0},
\end{equation}
for some constant $c_1$, $R$ large enough, and $\eta$ such that $\expecsub{\me^{8\eta\tauplus}}{\nu,0}<\infty$. For the last condition $8\eta<-\log(1-\pi_G)$ suffices, since the moment generating function of $\tauplus$ exists where that of $\tau$ exists; see the proof of Lemma~\ref{contraclem}. A similar last step as in that proof is taken, to conclude
\begin{equation}
\label{secmomalphahat}
\expecsub{(\hat\alpha-\alphastar)^2}{\nu,0} \leq \frac{2c_1}{R} \eta^2 \phi''(0).
\end{equation}

Now, in order to show that $\expec{\vnorm{\hat\mu-\mustar}^2} \leq c\,\vnorm{\nu-\mustar}^2$ with $c<1$ as required consider, with $\alpha$ fixed for the moment:
\begin{equation}
\begin{split}
\expecsub{\left (\hat\mu_x(\alpha)-\mu_x( \alphastar) \right)^2}{\nu,x}
& = [\mu_x(\alpha)-\mu_x( \alphastar)]^2 + \frac1R \varsub{Y_\nu(\alpha)}{\nu,x} \\
&\leq \vnorm{\mu(\alpha)-\mustar}^2 + \frac{M}{R} \vnorm{\nu-\mu(\alpha)}^2 \\
& \leq \left ( 1 + \frac{2M}{R} \right ) \vnorm{\mu(\alpha)-\mustar)}^2 + \frac{2M}{R} \vnorm{\nu-\mustar}^2,
\end{split}
\end{equation}
where the second step follows from Lemma~\ref{contraclem} (fixed $\alpha$), and the last step from the triangle inequality and Cauchy-Schwarz.
The Lipschitz property implies
\begin{equation}
 \expecsub{\left (\hat\mu_x(\alpha)-\mu_x( \alphastar) \right)^2}{\nu,x} \leq 
\left ( 1 + \frac{2M}{R} \right ) \partial^2 (\alpha-\alphastar)^2 + \frac{2M}{R} \vnorm{\nu-\mustar}^2.
\end{equation}
Now inserting $\alpha=\hat\alpha$ and taking expectations (noting that $\hat\mu$ and $\hat\alpha$ are independent), applying~\eqref{secmomalphahat}, one obtains 
\begin{equation}
\label{MSEmuhat}
\expecsub{\left (\hat\mu_x(\hat\alpha)-\mu_x( \alphastar) \right)^2}{\nu,x} \leq \frac{c_3}{R} \vnorm{\nu-\mustar}^2,
\end{equation}
for some constant $c_3$. Hence, if $R$ is large enough, the contraction condition~\ref{ass1} is fulfilled.

\medskip

As before, two timesteps are considered for condition~\ref{ass4}. Assume $\mu^{(0)}=\nu$, with $\nu$ arbitrary. Since 
$\hat\mu_x(\alpha) \leq \hat\mu_x(\alphamax)$ for $0\leq \alpha \leq \alphamax$ and any $x$, by the reasoning as in the proof of Lemma~\ref{twosteps}, 
an integer $R_1$ exist such that, if the number of replications is $R_1$, then 
\begin{equation}
\label{Markov}
\prc{\mu^{(1)} \in B}{\mu^{(0)}=\nu} \geq \frac12\quad\text{a.s.,}
\end{equation}
for $B$ 
as defined in Section~\ref{bounds} (note that the value of $\hat\alpha$ does not matter).

Fixing $\alpha=\alphamax$, it follows from Lemma~\ref{UIoverB} that 
$\{{\cal L}(\text{$Y_{0,\nu}(\alphamax)$ under ${\mathbb Q}_\nu$}): \nu \in B\}$ is uniformly integrable; this 
remains true if $\alphamax$ is replaced by a smaller $\alpha$. Now note, since $\hat\mu_0'(\xi) \geq \hat\mu_0'(0) 
\geq \hat\mu_0(0)$, that \eqref{sqerr} implies
\begin{equation}
|\hat\alpha-\alphastar| \leq \frac{|\hat\mu_0(\alphastar)-1|}{\hat\mu_0(0)}.
\end{equation}
As before, the ULLN applies, whence for all $\epsilon_1>0$ a sample size $R_2$ exists that guarantees
$\prc{|\hat\alpha-\alphastar|< \epsilon_1}{\mu^{(1)}}\geq 3/4$ a.s., uniformly on $\{ \mu^{(1)} \in B \}$.

Since $\{{\cal L}(\text{$Y_{0,\nu}(\alpha)$ under ${\mathbb Q}_\nu$}): \nu \in B, 0 \leq \alpha \leq \alphamax\}$ is uniformly integrable,
the ULLN implies for any $\epsilon_2>0$ the existence of a sample size $R_3$ such that
$\prc{|\hat\mu_x(\alpha) - \mu_x(\alpha)| < \epsilon_2} 
{\hat\alpha=\alpha, \mu^{(1)}} \geq 3/4$, uniformly on $\{ \mu^{(1)} \in B, 0 \leq \hat\alpha \leq \alphamax \}$.

Combining the last two results: for any $\epsilon_1, \epsilon_2>0$ a sample size $R_4$ exists such that
$\prc{|\hat\mu_x(\hat\alpha) - \mu_x(\hat\alpha)| < \epsilon_2, 1 \leq x \leq d, 
|\hat\alpha-\alphastar|<\epsilon_1}{\mu^{(1)}} \geq 1/2$, uniformly on $\{ \mu^{(1)} \in B \}$.
On the event in the last probability:
\begin{equation}
 \vnorm{\mu^{(2)} - \mustar} \leq \vnorm{\mu^{(2)} - \mu(\hat\alpha)} + \vnorm{\mu(\hat\alpha) - \mu(\alphastar)} < 
\epsilon_2+\partial\epsilon_1.
\end{equation}
Hence, if one chooses $\epsilon_1=\epsilon/2\partial$ and $\epsilon_2=\epsilon/2$, then
$\prc{\vnorm{\mu^{(2)}-\mustar}<\epsilon}{\mu^{(1)}} \geq 1/2$, uniformly on $\{ \mu^{(1)} \in B \}$.
Finally, combining this with~\eqref{Markov}, this proves that for any $\epsilon>0$ a sample size $R_0$ exists such that
$\prc{\vnorm{\mu^{(2)}-\mustar}<\epsilon}{\mu^{(0)}=\nu} \geq 1/4$, for any $\nu$.~\qed


\section{The general state space Markov chain perspective}
\label{MCperspective}
Desai and Glynn~\cite[\S 4]{desai2001simulation} discuss the convergence of the algorithms in \cite{KBCP99,baggerly2000exponential,desai2001adaptive} from the perspective of general
state space Markov chain theory, noting in particular the apparent impossibility to provide a proof based on this theory. They show that the Markov chains involved are not even Harris chains by proving that convergence in total variation to a Dirac measure at the zero-variance point does in fact \emph{not} happen. 
As an alternative method, they discuss the merits of Lyapunov approaches, hypothesizing that it might suffice to show that every $\epsilon$-neighborhood of the zero variance point is visited infinitely often, combined with weak continuity of the Markov kernel governing the process of iterates~$\mu^{(n)}$. As a negative answer they present a simple example that disproves their hypothesis. 
It is a Markov chain $X=\{X_n: n \geq 0 \}$ on $[0,\infty)$ with transition probabilities $P(x,1)=p(x)$, $P(x,x/2)=1-p(x)$ for $x>0$ and $P(0,0)=1$, 
where $p$ is a continuous function with $0 \leq p(x) \leq 1$ and $p(0)=0$. The kernel is weakly continuous and the Dirac measure at $0$ is a stationary distribution of this chain.  
They show that $X$ visits $(0,\epsilon)$ infinitely often, for every $\epsilon>0$ and every starting point $x \in (0,\infty)$. 
For $p$ with $\sum_{n=0}^\infty p(2^{-n})=\infty$, e.g., $p(2^{-j})=1/j$, the chain visits state 1 infinitely often, whence $X$ fails to converge to $0$.

Looking at the example from the perspective of Theorem~\ref{expconMC}, it can be proved that under the condition on $p$ no function $g$ exists that satisfies condition~\ref{ass1}: 
the contraction property forces $g(x)<0$ for $x$ near zero. Under some mild restrictions on $p$ (sufficient are: $\sup_{x \leq 1}p(x)<1$, $\inf_{x \geq 2} p(x)>0$), 
the other conditions (in particular~\ref{ass4}) are satisfied. With respect to the condition on $p$ it is close to an if and only if: 
requiring $\sum_{n=0}^\infty p(2^{-n})(1+\epsilon)^n<\infty$ for some $\epsilon>0$ is sufficient for the existence of a $g$ satisfying condition~\ref{ass1}.
This seems to indicate that the crucial property that made the counterexample fail to converge, is the drift away from zero that occurs near zero.

Theorem~\ref{expconMC} gives an alternate structure for a general convergence proof, consisting of a contraction property valid when the process has entered a ``good'' set $G$
and a uniform positive probability of entering $F \subset G$ in two steps. Part of the second step involves weak continuity: Lemma~\ref{weakcont} establishes weak continuity of the collection ($Y_\nu$ under ${\mathbb Q}_\nu$) for $\nu$ in a compact set $B$; this implies the weak continuity of the Markov kernel that governs the updating process, on the same set.
On the other hand, this weak continuity in itself is insufficient for this second step, the uniform integrability seems indispensible.
\bibliographystyle{plain}
\bibliography{zerovariance,bibAB}
\end{document}

%% file: wi4202mac.tex
\newcommand{\pr}{\ensuremath{{\mathbb P}}}
\newcommand{\prob}[1]{\ensuremath{{\mathbb P}\!\left( #1 \right)}}
\newcommand{\prc}[2]{\ensuremath{{\mathbb P}( #1 \, |\, #2)}}
\newcommand{\expec}[1]{\ensuremath{{\mathbb E}\mspace{-1mu}\left[#1\right]}}

\newcommand{\realR}{\ensuremath{\mathbb{R}}}
\newcommand{\me}{\mathrm{e}}
\newcommand{\md}{\mathrm{d}}

\newcommand{\probsub}[2]{\ensuremath{{\mathbb P}_{#2}\!\left( #1 \right)}}
\newcommand{\expecsub}[2]{\ensuremath{{\mathbb E}_{#2}\mspace{-1mu}\left[#1\right]}}

\newcommand{\cexpec}[2]{\ensuremath{{\mathbb E}\left [ #1 \, |\, #2 \right ]}}

\newcommand{\eqdist}{\stackrel{d}{=}}

\newcommand{\convas}{\stackrel{\rm a.s.}{\longrightarrow}}

\newcommand{\seq}[3]{\ensuremath{{#1}_{#2},\ldots, {#1}_{#3}}}  
\newcommand{\infseq}[3]{\ensuremath{{#1}_{#2}, {#1}_{#3}, \ldots}}  

\newcommand{\indicator}[1]{\ensuremath{{\bf 1}_{\{#1\}}}}
\newcommand{\vnorm}[1]{\ensuremath{|| #1 ||}}

